\documentclass[11pt,sts,noinfoline]{imsart}

\RequirePackage[OT1]{fontenc}
\RequirePackage{amsthm,amsmath,amssymb}
\usepackage[numbers,sort&compress]{natbib}
\RequirePackage{hyperref}
\usepackage{float}
\usepackage{mathrsfs}
\usepackage{verbatim}
\usepackage{graphics}	
\usepackage{dcolumn}	
\usepackage{bm}		
\usepackage{bbm,upgreek}
\usepackage{mathtools}
\usepackage{color}
\usepackage[cal=boondoxo]{mathalfa}
\usepackage{dsfont}
\usepackage{cancel}
\usepackage{xcolor}
\usepackage{harpoon}
\usepackage{accents}
\usepackage{tabularx,enumitem}

\arxiv{math.PR/00000000}

\startlocaldefs
\numberwithin{equation}{section}
\theoremstyle{plain}

\endlocaldefs

\newcommand{\paren}[1]{\left( #1 \right)}
\newcommand{\integral}[2]{\underset{#2}{\int} \mathrm{d} #1 \,} 

\newcommand{\expectation}[2]{\underset{#2}{\mathbb{E}_{#1}} \,}

\newcommand{\param}[0]{\bm \theta}

\newcommand{\MLE}{\hat{\bm \theta}}

\newcommand{\prior}{\pi}

\usepackage{array}
\newcolumntype{L}[1]{>{\raggedright\let\newline\\\arraybackslash\hspace{0pt}}m{#1}}
\newcolumntype{C}[1]{>{\centering\let\newline\\\arraybackslash\hspace{0pt}}m{#1}}
\newcolumntype{R}[1]{>{\raggedleft\let\newline\\\arraybackslash\hspace{0pt}}m{#1}}

\newcommand*{\inliner}[1]{{\smash{\Scale[0.85]{#1}}}}

\newcommand*{\Scale}[2][4]{\scalebox{#1}{$#2$}}%

\newcommand{\idea}[1]{\subsection{ #1}}

\begin{document}

\begin{frontmatter}
\title{The Lindley paradox: the loss of resolution in Bayesian inference.}
\runtitle{Optimal resolution in Bayesian inference}

\begin{aug}
\author{\fnms{Colin H.} \snm{LaMont}\ead[label=e1]{lamontc@uw.edu}}
\and
\author{\fnms{Paul A.} \snm{Wiggins}\ead[label=e2]{pwiggins@uw.edu}}

\address{University of Washington\\
Seattle, WA\\
\printead{e1,e2}}

\runauthor{C. H. LaMont and P. A. Wiggins}

\affiliation{University of Washington}

\end{aug}

\begin{abstract}
There are three principle paradigms of statistics: Bayesian, frequentist and information-based inference. Although these paradigms are in agreement in some contexts, the Lindley paradox describes a class of problems, models of unknown dimension, where conflicting conclusions are generated by frequentist and Bayesian inference. This conflict can  materially affect the scientific conclusions. Understanding the Lindley paradox---where it applies, why it occurs, and how it can be avoided---is therefore essential to the understanding of statistical analysis. In this paper, we revisit the Lindley paradox in the context of a simple biophysical application. We describe how predictive and postdictive measures of model performance provide a natural framework for understanding the Lindley paradox. We then identify methods which result in optimal experimental resolution for discovery.
\end{abstract}

\footnote{AMS classification: 62F15, Keywords: Lindley paradox; information criterion; Bayesian inference; pseudo-Bayes Factors}

\end{frontmatter}

\section{Introduction}

With advances in computing, Bayesian methods have experienced a strong resurgence. 
Proponents of Bayesian inference cite numerous practical and philosophical advantages of the paradigm over classical (frequentist) statistics \cite{Bernardo1994}.
The most compelling argument in favor of Bayesian methods is the natural hedging between competing hypotheses and parameter values. This hedging mechanism (i.e. model averaging) protects against over-fitting in singular models and has led to excellent performance in machine learning applications and many other contexts, especially those which require the synthesis of many forms of evidence \cite{watanabe2009,Bernardo1994}. But the practical and philosophical problems that motivated the development of frequentist methods remain unresolved: (i) There is  no commonly agreed upon procedure for specifying the Bayesian prior  and (ii) statistical inference can depend strongly upon the prior. This dependence creates a discrepancy between Bayesian and frequentist methods: the Lindley paradox.


We analyze  Bayesian model selection with respect to the relative partition of information between the data and the prior. This analysis leads to novel connections between Bayesian, information-based, and frequentist methods.  We demonstrate that a large prior information partition results in model selection consistent with the Akaike Information Criterion (AIC) \cite{akaike1773}, while the opposite limit of the information partition results in model selection consistent with the Bayesian Information Criterion (BIC) \cite{Schwarz1978a}. Intermediate partitions interpolate between these well known limits. Although the AIC limit is well defined and robust, the BIC limit depends sensitively on the \textit{ad hoc} definition of a \textit{single measurement}. Furthermore, the BIC limit corresponds to a loss of resolution. This loss of resolution might result in the unnecessary purchase of more sensitive equipment or the collection of unreasonable sample sizes. 

As a result, we question the suitability of BIC model selection (or Bayesian inference with an uninformative prior) at finite sample size.
The large-prior-information  regime of Bayesian inference can be achieved in almost any practical Bayesian implementation by the use of pseudo-Bayes factors\cite{geisser1979,gelfand1992}.  This approach circumvents the Lindley paradox while maintaining many advantages of Bayesian inference.




\begin{figure}
\includegraphics[width=.7\textwidth]{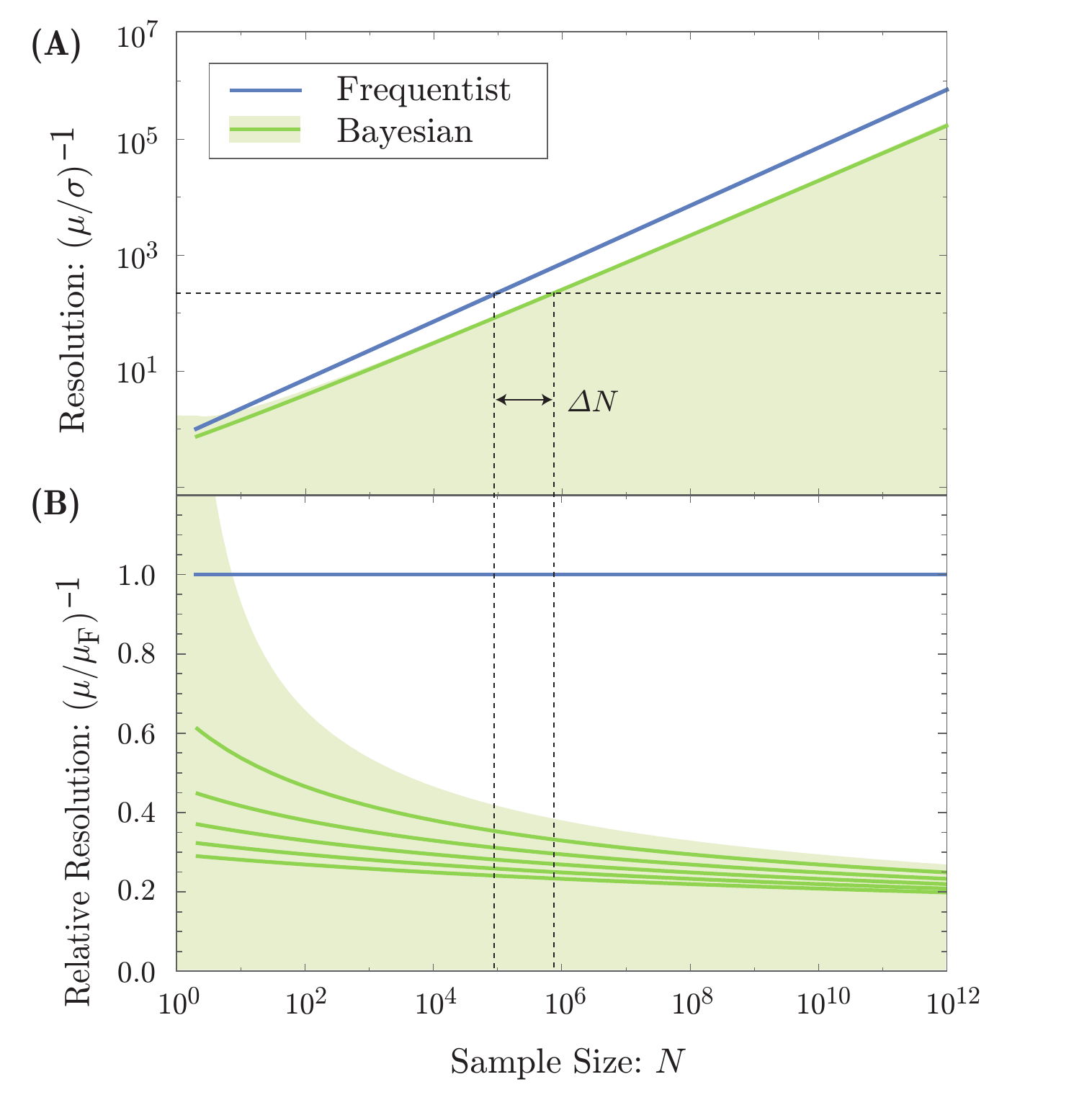}
\caption{\textbf{Loss of resolution in Bayesian inference.} \textbf{Panel A:} The resolution on detected bead displacement (the alternative hypothesis) is plotted as a function of sample size $N$. The increase in resolution is due to the decrease in the error in the mean $\sigma_\mu=\sigma/\sqrt{N}$. The resolution of both frequentist and Bayesian inference increase, but the frequentist resolution is higher. A dotted line represents the size of a putative displacement. The frequentist analysis detects this displacement at a smaller sample size than the Bayesian analysis.
\textbf{Panel B:} To draw attention to the difference between Bayesian and frequentist resolution, we plot the resolution relative to the frequentist resolution $\mu_F$. To illustrate the prior dependence of the Bayesian analysis, we have drawn curves  corresponding to various choices of prior volume $V_0$. 
\label{fig:Res}}
\end{figure}

\subsection{A simple example of the Lindley paradox}
A simple example emphasizes the difference between  Bayesian and frequentist forms of statistical support for models of unknown dimension. 
Suppose an observer measures the position of a bead in the course of a biophysics experiment.
The position is first determined with negligible uncertainty. After a perturbation is applied, $N$ measurements are made of the bead position: $\inliner{ x^N \equiv (x_1,x_2,...,x_N)}$.
The $N$ measurements are assumed to be independent and identically distributed (iid) in a normal distribution centered on the unknown true displacement $\mu$ with known variance $\sigma^2$ where $\mu=0$ if the bead is unmoved and $\mu \ne 0$ otherwise. 

In the Bayesian paradigm, we  must specify priors over the parameters for the two models $\pi_i({\bm \theta})$. Model zero (null hypothesis) is parameter free since $\mu=0$, but  model one (alternative hypothesis) is parameterized by unknown mean $\mu$. The true value $\mu_0$ is unknown and to represent this ignorance, we use a \textit{vague} conjugate prior, choosing a normal prior centered on zero with a  large variance $\tau^2$. A canonical objective Bayesian approach to model selection is to assume the competing models have equal prior probability. The model with the largest  posterior probability is selected. The  experimental resolution for detecting a change in the bead position is then: 
\begin{equation}
|\hat{\mu}| >\sigma_\mu \cdot \sqrt{2\log \tau/\sigma_\mu}, \label{eqn:bayesres}
\end{equation}
while the frequentist \textit{rule of thumb} ($\approx 95\%$ confidence level) for rejecting the null hypothesis is:
\begin{equation}
|\hat{\mu}| > \sigma_\mu \cdot 2, \label{eqn:freqres}
\end{equation} 
where $\sigma_\mu \equiv \sigma/\sqrt{N}$ is the uncertainty in $\mu$.
The difference between the conditions defined by Eqns.~\ref{eqn:bayesres} and \ref{eqn:freqres}
reveals that the paradigms may come to conflicting conclusions about model selection, as illustrated in Fig.~\ref{fig:Res}. D.~Lindley emphasized this conflict by describing the following scenario: If the alternative hypothesis is true, for a suitable choice of $\tau$ and  sample size $N$, the null hypothesis could be simultaneously (i) rejected at a 95\% confidence level \textit{and} (ii) have 95\% posterior probability  \cite{Lindley1957}!  This conflict between statistical paradigms has been called the Lindley paradox.

Many practitioners of Bayesian inference believe that priors may be a formal necessity but have minimal influence on inference. 
For instance, the posterior probability for $\mu$ is independent of the prior in the uninformative limit $\tau \rightarrow \infty$.
However, as we see in Eqn.~\ref{eqn:bayesres}, inference on model identity remains critically dependent on the prior (value of $\tau$). In the limit that $\tau \rightarrow \infty$, no finite observed displacement $\hat{\mu}$ is sufficient to support the alternative hypothesis that the bead has moved! This paradoxical condition is called the Bartlett paradox \cite{Bartlett1957}.

\section{Data partition}

\idea{The definition of frequentism and Bayesian paradigms} 
We wish to study a generalized class of decision rules that include methods from all three paradigms of inference. In the current context, we will use the log likelihood ratio: 
\begin{equation}
\lambda(x^N) \equiv h_0( x^N| \hat{\bm \theta}_x ) - h_1( x^N| \hat{\bm \theta}_x ), \label{eqn:teststat}
\end{equation}
as a \textit{frequentist test statistic} where $h$ is the Shannon information $h \equiv -\log q$ and $\hat{\bm \theta}_x$ is the  maximum likelihood estimate of the parameters of the respective model.  We shall define a decision rule:  
\begin{equation}
\lambda(x^N) < \lambda_*, \label{eqn:dr}
\end{equation}
to select model zero where $\lambda_*$ is the critical value of the test statistic. We will refer to the decision rule as \textit{frequentist} if $\lambda_*$ is sample-size independent in the large-sample-size limit of a regular model. This definition includes both the frequentist Neyman-Pearson likelihood ratio test as well as the information-based paradigm (AIC). In the Bayesian paradigm, we will define the decision rule in terms of the log-Bayes factor:
\begin{equation}
\lambda_B(x^N) \equiv h_0(x^N)-h_1(x^N),
\end{equation}
where $q(x^N)$ is the marginal likelihood and $h(x^N)$ is the respective Shannon information.  We define the decision rule: $\lambda_B(x^N)<0$ to select model zero. Although the Bayes factor is not  a test statistic---an orthodox Bayesian approach is to compute a posterior on model identity---the decision rule captures how the Bayes factor is typically used in practice.

In the large-sample-size limit, the Bayesian decision rule is equivalent to Eqn.~\ref{eqn:dr} with $\lambda_*$ proportional to $\log N$ to leading order in $N$. Therefore, we will define a decision rule Eqn.~\ref{eqn:dr} as \textit{Bayesian} if the critical test statistic $\lambda_*$ is sample-size dependent. This definition includes standard Bayesian model selection as well as the Bayesian information criterion (BIC).

\subsection{Prior information content}

The paradoxically-large displacement needed to select the alternative hypothesis is a consequence of the uninformative prior ($\tau \rightarrow \infty$). To be more precise about the descriptors  \textit{informative} and \textit{uninformative}, we can compute  expected-parameter-information content of the data set $x^{N}$ \cite{Lindley1956}:
\begin{equation}
I(x^{N}) \equiv \expectation{\bm \theta}{\pi({\cdot}|x^{N})}{\log \pi({\bm \theta}|x^{N})/\pi({\bm \theta}}), \label{eqn:prior_info_content}
\end{equation}
which is equal to the KL Divergence of the posterior and prior. $I\ge 0$ and will increase with sample size. Given $N$ new measurements, we call the prior $\pi$ \textit{uninformative} if $I$ is large. 

A standard approach to specify an informative prior is the elicitation of a prior from an expert \cite{Bernardo1994}. It is convenient to make the concrete assumption that the expert knowledge is the result of previous measurements, which we can write explicitly $x^{N_T}$. Our posterior on these measurements $\inliner{\pi({\bm \theta}|x^{N_T})}$ is computed from some suitably flat prior $\pi({\bm \theta})$.  The $x^{N_T}$ is then used to construct a new informative prior: 
\begin{equation}
\pi'({\bm \theta}) \equiv  \pi({\bm \theta}|x^{N_T}), \label{eqn:info}
\end{equation}
%
where the primed distributions are computed with respect to the informative prior (Eqn.~\ref{eqn:info}). This Bayesian update rule was concisely summarized by D.~Lindley: \textit{Today's posterior is tomorrow's prior.}  

Let the new measurements made be re-labeled $x^{N_G}$. We can re-compute the marginal likelihood $q'$ using the new prior $\pi'$. $q'$ has a second interpretation, the Bayesian predictive distribution computed from the original prior $\pi$: 
\begin{equation}
q'(x^{N_G}) = q(x^{N_G}|x^{N_T}) \equiv q(x^{N})/q(x^{N_T}), \label{eqn:evNG}
\end{equation}
where $x^{N}$ represents the entire data set of $N = N_G+N_T$ measurements.
This distribution is predictive since it predicts or \textit{generalizes} on data set $x^{N_G}$ given a \textit{training} data set $x^{N_T}$. Adjustment of the data partition between the training set (size $N_T$) and the generalization set (size $N_G$) can be understood as adjusting the information content of the prior. If $N_G \gg N_T$, the prior is uninformative relative to the data.

\idea{The Bayesian cross entropy} 
The general problem of predicting observations $x^{N_G}$, conditioned on $x^{N_T}$ where $N=N_G+N_T$ is closely related to a natural metric of performance: a predictive \textit{cross entropy} \cite{BurnhamBook}
\begin{equation}
H^{{N_G}|{N_T}} \equiv {\textstyle \frac{N}{N_G}}\ \inliner{\expectation{X}{p(\cdot)}} h(X^{N_G}|X^{N_{T}}), \label{eqn:cross}
\end{equation}
where $p(x^N)$ is the true distribution of observations $x^N$. The cross entropy is rescaled to correspond to the total sample size $N$.
We can view model inference using the evidence Eq.~\ref{eqn:evNG} as choosing the model which is estimated to have the optimal performance under this metric.
Since $H$ can only be computed if the true distribution $p$ is known, it will be useful to empirically estimate it.
A natural  estimator is the leave-k-out estimator \cite{gelfand1994}
\begin{equation}
\hat{H}^{N_G|N_T}(x^N) \equiv {\textstyle \frac{N }{N_G}}  {\expectation{X}{P\{x^N\}}} \, h(X^{N_G}|X^{N_T}) \label{eq:psbfestimator},
\end{equation}
where $\hat{H}$ estimates $H$ and the empirical expectation $\mathbb{E}$ is taken over all unique permutations of the observed data between the training and generalization sets.

This estimator uses \textit{cross validation}: there is no double use of data since the same observations never appear in both the generalization and training sets. Methods like empirical Bayes \cite{MacKay1992a,MacKay1992b,Aitkin1991}, where the prior is fit to the data to maximize the evidence, implicitly use the data twice and are therefore subject to the same over-fitting phenomenon as maximum likelihood estimation.

\idea{Pseudo-Bayes factors} The natural strategy would be to  compute the  model posterior probability (or Bayes factor) using the evidence $q'(x^{N_G})$. But, for small $N_G$, $h(x^{N_G}|x^{N_T})$ typically exhibits large statistical fluctuations since only a small fraction of the data $x^{N_G}$ is used for inference on the model identity even though there is more non-redundant information encoded in $x^{N_T}$. To re-capture this missing information, we replace $\inliner{h(x^{N_G}|x^{N_T})}$ with $\inliner{\hat{H}^{N_G|N_T}}$. Therefore, in analogy with the log Bayes factor, the log-{\em pseudo-Bayes factor} is defined \cite{gelfand1994}:
\begin{equation}
\lambda_{\rm PB}^{N_G|N_T}(x^N) \equiv \hat{H}_0^{N_G|N_T} - \hat{H}_1^{N_G|N_T}, \label{eq:psbf}
\end{equation}
which depends on the data partition $N_G$ and $N_T$. We define the decision rule: $\lambda_{PB}(x^N)<0$ to select model zero.

Two data partitions have been discussed in some detail. 
A maximal-traning-set limit, where $N_T = N-1$ and $N_G = 1$, corresponds to  Leave-one-out cross validation (LOOCV) and has been  studied extensively \cite{gelfand1994,geisser1979,gelfand1992,vehtari2002bayesian,Spiegelhalter2002}.
A mininimal-training-set limit has also been explored in which $N_T$ is as small as possible such that $\pi'$ is proper \cite{atkinson1978posterior,Smith1980}.

We focus on the example of a pairwise model selection to compare with canonical frequentist inference, but a selection among any number of models can be performed by selecting the model with the smallest cross-entropy estimator.

\subsection{Information Criteria} 
\begin{figure}
\includegraphics[width = .7\textwidth]{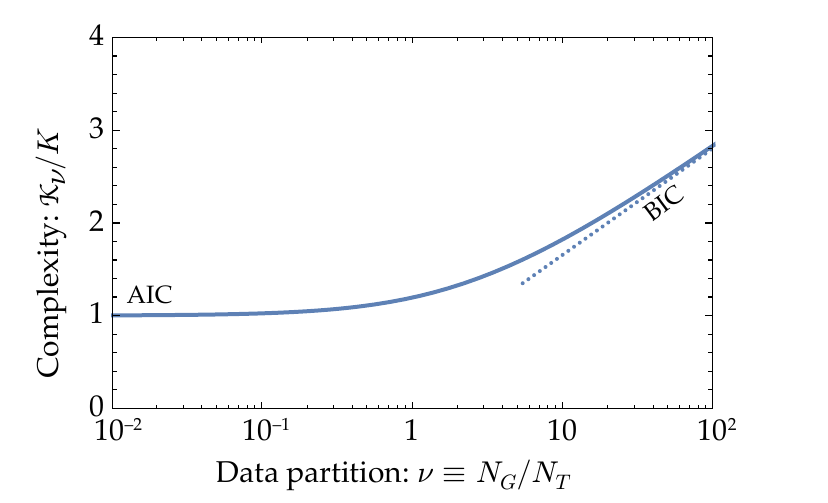}
\caption{{\bf Complexity as a function of data partition.} Complexity can  be understood as a penalty for model dimension $K$. The data partition parameter controls the relative amount of information in the prior. In a predictive limit ($\nu\rightarrow 0$), the training set is large compared with the generalization set and the complexity is small.
This is the AIC limit. At the other extreme ($\nu\rightarrow \infty$), all the data is partitioned into the generalization set and therefore the prior is uninformative and the complexity is large. This is the BIC limit. \label{valvtrain}}
\end{figure}

To systematically investigate the dependence of inference on the data partition in the pseudo-Bayes Factor, we propose a novel estimator of the cross entropy $H$ whose dependence on the data partition is explicit. The data partition will be parameterized by $\nu \equiv N_G/N_T.$
We define a generalized Information Criterion:
\begin{equation}
{\rm IC}^{\nu}(x^N) \equiv h(x^{N}|\MLE_x) + {\cal K}_{\nu},
\end{equation} 
where  the complexity ${\cal K}_{\nu}$ is the bias, chosen to make ${\rm IC}^{\nu}$ an unbiased estimator of $H^{N_G|N_T}$. The log-pseudo-Bayes Factor can be constructed using the information criterion. The information criterion is typically much easier to evaluate than the leave-k-out formulation. Since the first term in the definition of ${\rm IC}^{\nu}$ is independent of $\nu$, the data-partition dependence is completely characterized by the complexity ${\cal K}_{\nu}$.

Assuming $\pi$ is uninformative and $\inliner{q(x^N|{\bm \theta})}$ is a regular model in the large sample size limit, the Laplace approximation holds and the complexity has a simple form:
\begin{equation}
{\cal K}_\nu = {\textstyle \frac{1}{2}} K \left[1+(1+\nu^{-1}) \log (1 + \nu)\right] \label{compnu},
\end{equation}
which is only a function of the parameter-space dimension $K$ and the data partition $\nu$.
The complexity is plotted as a function of the data partition $\nu$ in Fig.~\ref{valvtrain}.

\subsection{Decision rules and resolution}
\begin{figure}
\includegraphics[width=.7\textwidth]{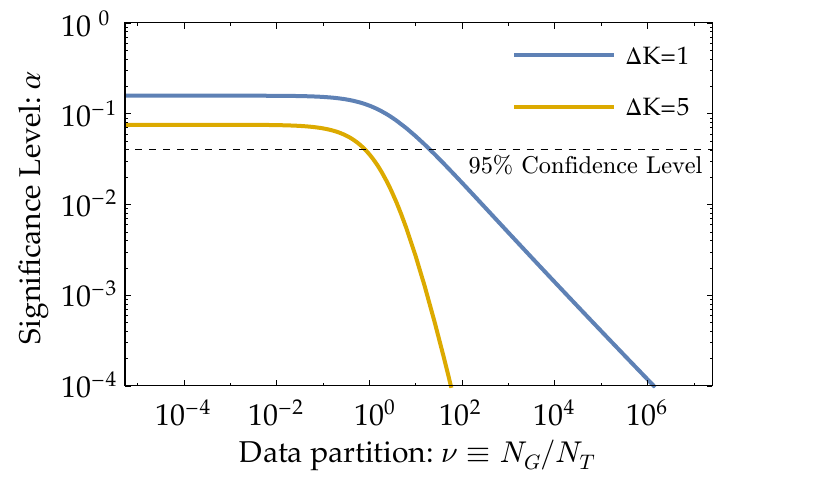} 
\caption{{\bf Significance level as a function of data partition.} To make an explicit connection between the frequentist significance level and the data partition, it is useful to compute the significance level implied by the predictive decision rule.  
In the postdictive regime, corresponding to an uninformative prior, the significance level steeply decreases with increasing $\nu$, resulting in a strong Lindley paradox. \label{fig:pval}}
\end{figure}

With the information criterion above, we can connect a (pseudo-)Bayes factor with a particular data partition $\nu$ to an effective decision rule. We choose model one if
\begin{equation}
h_0(x^{N}|\MLE_x) - h_1(x^{N}|\MLE_x) >\Delta {\cal K }_\nu,
\end{equation}
where $\Delta {\cal K }_\nu$ is the difference in the complexity of the models.  We can also connect these decision rules to choices of a frequentist significance level as described in the supplement. A plot of this function for two different values of $\Delta K$ is shown in Fig.~\ref{fig:pval}.
Of particular practical experimental importance is the minimal signal to noise ratio at which our decision rule will choose a larger model. Returning to the biophysical problem described in the introduction, the minimal resolvable bead displacement is 
\begin{equation}
|\hat{\mu}| > \sigma_\mu\cdot\sqrt{1+(1+\nu^{-1}) \log (1 + \nu)}
\label{eqn:resolution}
\end{equation} 
where the RHS is the inverse resolution. The resolution is monotonically decreasing  in $\nu$. The smallest $\nu$ gives us the highest resolution.

\section{The Lindley paradox}

The Lindley paradox can be understood from the perspective of the relative partition of information between the prior and the data. By defining the complexity (Eqn.~\ref{compnu}), we can explore the partitioning of this data by studying the decision rule and resolution as a function of the partition $\nu$.

\subsection{Classical Bayes and the Bartlett paradox}

\begin{figure}[t]
\includegraphics[width=.5\textwidth]{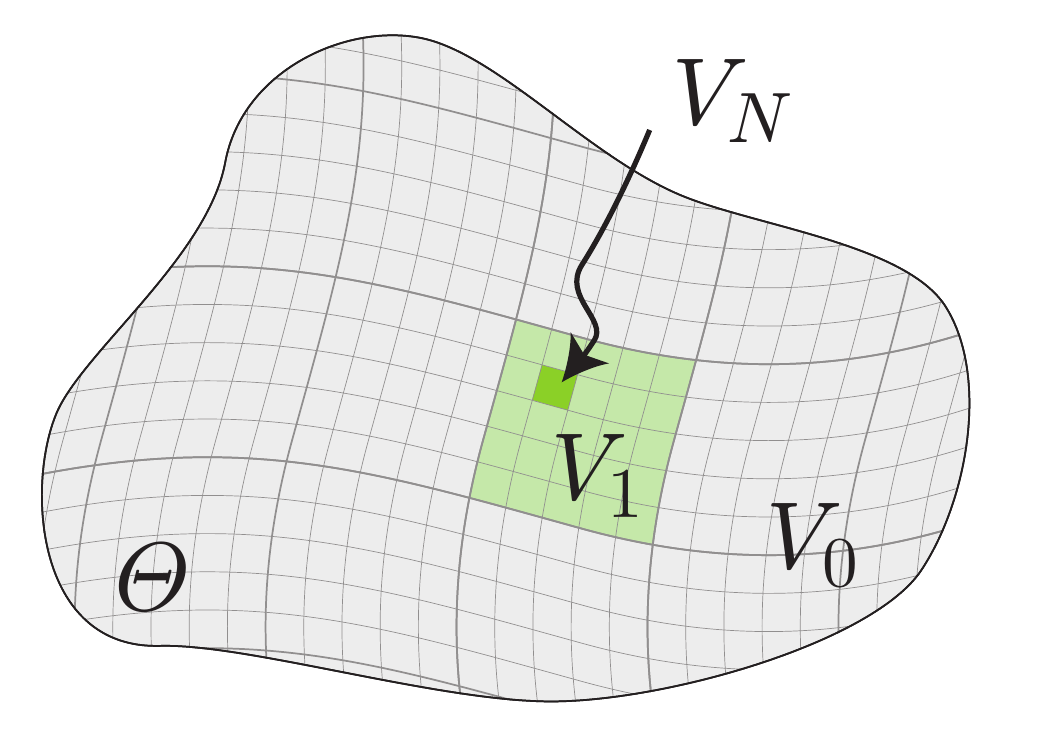}
\caption{\textbf{The geometry of the Occam factor.} The total volume of plausible parameter values for a model is $V_0$. The volume of acceptable parameter values after a single measurement is $V_1$. The volume of acceptable parameter values after  $N$ measurements is $V_N$. The Occam factor is defined as the probability of randomly drawing a parameter from initial volume $V_0$ consistent with the $N$ measurements: ${\rm Pr} \approx {\cal N}^{-1}$ where ${\cal N}\equiv V_0/V_N$ is the number of distinguishable distributions after $N$ measurements. Lower dimension models are naturally favored by the Occam factor since the number of distinguishable models $\cal N$ is smaller.  \label{fig:volume}} 
\end{figure}
For the classical Bayes factor, $N_T =0$ or $\nu\rightarrow \infty$.  If the prior is flat on an infinite-volume parameter manifold, the complexity  ${\cal K}_\nu$ becomes infinite. This scenario always favors the smaller model, regardless of the goodness of fit, resulting in the Bartlett paradox. 

If the parameter-manifold volume $V_0$ is finite, so is the complexity. In the large sample size limit, the marginal likelihood can be written in an intuitive form: 
\begin{equation}
q(x^N) = \textstyle \frac{V_N}{V_0} \times q(x^N|\MLE_x),
\end{equation} 
where $V_N$ is the volume of parameter manifold consistent with the data $x^N$ and $\inliner{\MLE_x}$ is the maximum likelihood estimate of the parameters. (We define this volume more precisely in the supplement.)
The first factor on the RHS is the Occam factor or the probability of randomly drawing a parameter (consistent with $x^N$) from the prior distribution $\pi$. 
Complex models (large $K$) with uninformative priors have \textit{small} Occam Factors, due to the large volume of plausible parameters ($V_0$), relative to the volume of the parameter manifold consistent with the observations ($V_N$). Large Occam factors give rise to a natural mathematical realization of the Occam Razor: \textit{Among competing hypotheses, the one with the fewest assumptions [parameters] should be selected} \cite{MacKay1992a}. This effect is illustrated schematically in Fig.~\ref{fig:volume}. 
Both infinite and finite-but-large-volume parameter manifolds can give rise to strong Lindley paradoxes.


\subsection{Minimal training set and Lindley Paradox}

We might use a \textit{minimal training set} to remove the dependence on the potentially divergent volume $V_0$  \cite{Berger1996} which corresponds to the large-data-partition limit: $\nu\gg1$.
It is difficult to define this minimal training set in a satisfactory way \cite{ohagan1995}.
The most natural option is to set $N_T=1$ and $N_G=N-1$, which results in the Bayesian information criterion (BIC) \cite{Schwarz1978a}:
\begin{equation}
{\rm BIC}(x^N) \equiv h(x^{N}|\MLE_x) + {\textstyle \frac{ K}{2}} \log N,
\end{equation}
in the large-sample-size limit. We can compute a limit on the smallest resolvable change in position:
\begin{equation}
\left|\hat{\mu}\right| > \sigma_\mu \cdot \sqrt{\log N}, \label{eqn:BICres}
\end{equation}
which is free from the \textit{ad hoc} volume $V_0$ of the uninformative prior.  This approach  resolves the Bartlett paradox, but leads to a strong Lindley paradox---conflict with Frequentist methods in some critical range of sample sizes.

The $\log N$ dependence of BIC results in some troubling properties. If we now bin pairs of data points, the empirical mean $\hat{\mu}$ and the standard error $\sigma_\mu$ are unchanged, but $N \rightarrow N/2$, changing the complexity and therefore the decision rule and resolution. 
Therefore, although BIC does not depend on the choice of prior support, it does depend on an \textit{ad hoc} choice as to what constitutes a single sample.

\subsection{Frequentist prescription and AIC} 
The complementary limit describes a maximal training set ($\nu \rightarrow 0$). 
In this limit, ${\rm IC}^\nu$  corresponds to the Akaike Information Criterion (AIC):
\begin{equation}
{\rm AIC}(x^N) \equiv h(x^{N}|\MLE_x) + K,
\end{equation}
where the complexity is equal to the dimension of the parameter manifold  $K$. This leads to a sample-size-independent critical value of the test statistic in Eqn.~\ref{eqn:dr}, and is therefore {\em frequentist}.
Like the log-Occam factor, $K$ can be reinterpreted as a penalty for model complexity that gives rise to a distinct information-based realization of the Occam Razor: \textit{Parsimony implies predictivity.}

The smallest resolvable change in position using AIC is
\begin{equation}
\left|\hat{\mu}\right| > \sigma_\mu \cdot \sqrt{2}. \label{eqn:aicres}
\end{equation} 
AIC can also be viewed as the performance of the model against the next observation $x_{N+1}$ \cite{BurnhamBook}. For this reason, $\nu\rightarrow 0$ is called the \textit{predictive limit}. The canonical Bayesian approach estimates the marginal likelihood of the observed data $x^N$ from the prior. It is therefore postdictive. We therefore call $\nu\rightarrow \infty$ the \textit{postdictive limit}.
The AIC and BIC penalties are often used as complimentary heuristics in model selection \cite{presse2016}.
Our cross entropy description shows that they can be interpreted as Bayesian objects which differ only in the choice of data partition $\nu$.

The predictive limit is expected to result in maximum resolution and consistent inference (\textit{i.e.} independent of the prior and data partition).  Unlike BIC, it is essentially independent of data binning in the large sample size limit. (See Fig.~\ref{valvtrain}.) Although AIC is computed using a point estimate, a pseudo-Bayes factor for  $N_G=1$ and $N_T=N-1$ (\textit{i.e.}~LOOCV) corresponds to the same predictive limit.

\begin{figure}
\includegraphics[width=.7\textwidth]{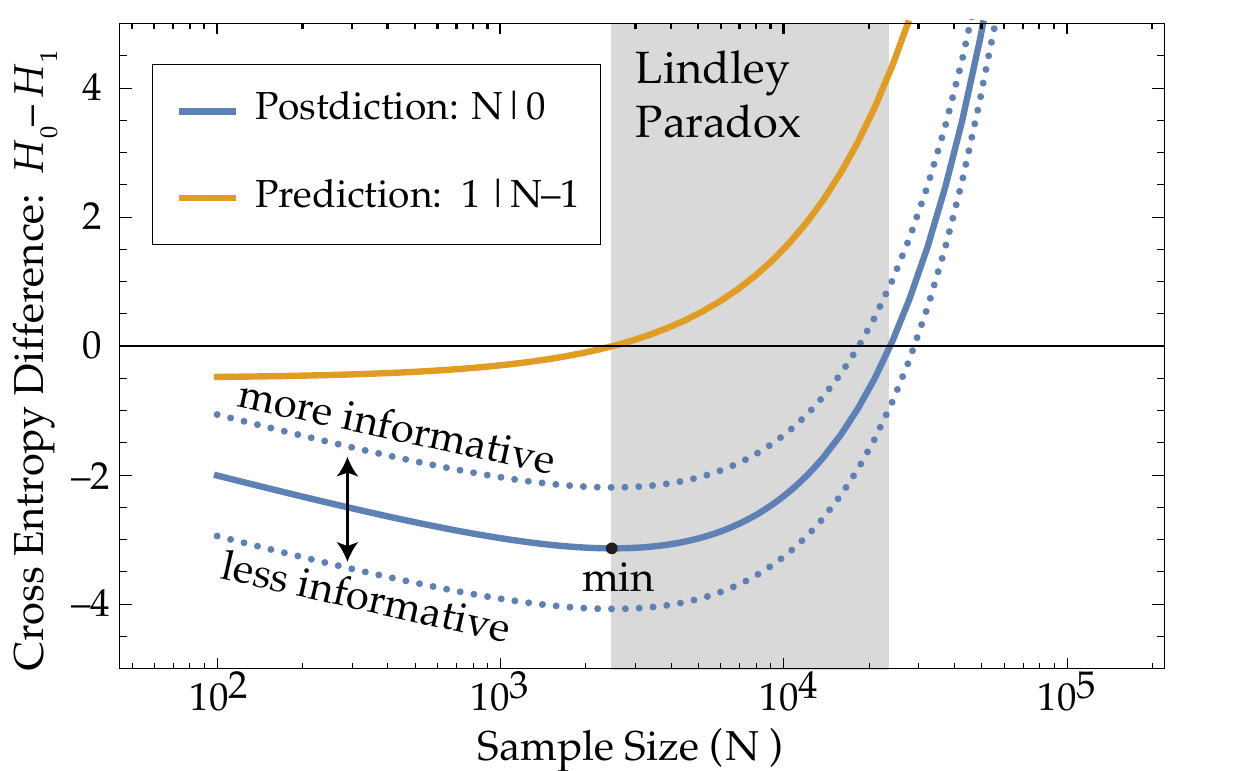} 
\caption{{\bf Visualizing the pre and postdictive decision rules.} The cross entropy difference for Bayesian inference (postdiction: $N|0$) and the predictive limit ($1|N-1$) are plotted as a function of sample size $N$. $\Delta H>0$ results in the selection of the alternative hypothesis. Both measures initially favor the null hypothesis. The use of a more (less) informative prior causes the postdictive curve to be shifted up (down). Since the predictive $H$ is the derivative of the postdictive $H$, the prior  does not influence the inference of the predictive observer.
The predictive  curve crosses zero first, leading the predictive observer to support the alternative hypothesis. Since the predictive $H$ is the derivative of the postdictive $H$ with respect to $N$, the sample size at which the predictive observer switches to the alternative hypothesis corresponds to the sample size where the postdictive observer has the most evidence for the null hypothesis. The two measures are in conflict (grey region) until the postdictive $H$ crosses zero at a significantly larger sample size $N$. The Bayesian therefore requires significantly more evidence to reach the same conclusion as the predictive observer.  
\label{Fig:Lindely}}
\end{figure}

\idea{Log evidence versus AIC} 
\label{sec:dd}

A convenient heuristic for understanding the relation between AIC and the log evidence can be understood from the relation between the cross entropy using $\inliner{H^{1|N-1}}$ and $\inliner{H^{N|0}}$ which are estimated by AIC and $h(x^N)$ respectively. If we approximate the finite difference in $\inliner{H^{1|N-1}}$ as a derivative in the large $N$ limit, $\inliner{H^{1|N-1}}$ can be approximated:
\begin{equation}
H^{1|N-1} = N \partial_N H^{N|0}  +{\cal O}(N^{-1}).
\end{equation}
Therefore, we can understand the relation between the conflicting information criteria AIC and BIC in the following way: AIC is the derivative of BIC.

This heuristic can naturally explain why AIC is free from the strong prior dependence which leads to the Lindley paradox. In the context of an uninformative prior, the expected log evidence $\inliner{H^{N|0}}$ has an ambiguous offset corresponding to the prior choice, leading different individuals to make different inference on the model identity. 
$\inliner{H^{1|N-1}}$, estimated by AIC, is independent of the unknown constant since the slope of $\inliner{H^{N|0}}$ is independent of its offset. This relationship is illustrated schematically in Fig.~\ref{Fig:Lindely}.

A second interesting feature of the heuristic relates to the sample size dependence of the predictive and postdictive decision rules.
The sample size at which the predictive statistician begins to favor the alternative hypothesis corresponds to the same sample size at which the postdictive statistician has maximum confidence in the null hypothesis!  (See Fig.~\ref{Fig:Lindely}.) The difference between a function and its derivative explains both the connection and inconsistency of the predictive and postdictive decision rules.


\subsection{When will a Lindley paradox occur?}
We stated  that the Lindley paradox is a consequence of an insufficiently informative prior, but we have studied  differences in the performance of predictive and postdictive  decision rules.  We now discuss the connection between these  equivalent formulations. Let us define the difference between the predictive and postdictive cross entropy:
\begin{equation}
I'(x^N) \equiv \hat{H}^{1|N-1}-\hat{H}^{N|0}. \label{eqn:missingInfo}
\end{equation}
In the large-sample-size limit, 
we can express $I'$ in terms of the expected-parameter-information content of the  data $I$: 
\begin{equation}
I' = I(x^{N}) + {\cal O}(N^0),
\end{equation}
as shown in the supplement. $I'$, the mismatch between pre and postdictive measures of performance,  can be interpreted as the \textit{missing information} from an uninformative  prior $I$. 
The missing information is only \textit{missing} before sample $x^N$ is observed. A model may be extremely predictive, even if the missing information was infinite, once $x^N$ has been observed.



\section{Discussion}

By defining a novel information criterion that estimates the cross-entropy, we established a continuous bridge between canonical Bayesian and information-based model selection defined in terms of a data-partition between training and generalization data sets. The strength of the Lindley paradox, the mismatch between Bayesian and frequentist inference on hypotheses, can be re-interpreted in terms of prior information content (\textit{i.e.}~the data partition). We studied the properties of model selection with respect to the data partition.
Two solutions to the Lindley paradox have been widely discussed: (i) adapt the frequentist paradigm by making the significance level sample-size dependent \cite{Berger2003} or (ii) adapt the Bayesian paradigm by making the prior sample-size dependent. We  advocate for taking the second approach.

\idea{A canonical Bayesian perspective on the Lindley paradox}

It is important to acknowledge that the Bayesian perspective on the Lindley paradox is valid. Returning to the  biophysical example, if we interpret the alternative hypothesis \textit{precisely}, we define a uniform prior probability density over an infinite-volume manifold in the uninformative limit ($\tau$ and $V_0 \rightarrow \infty$).  Therefore, the \textit{a priori} probability of picking a displacement consistent with the data ($\approx V_N/V_0)$ is vanishingly small in the alternative hypothesis.
\textit{Fine tuning} would be required to make $\hat{\mu}$ finite and therefore the null hypothesis is strongly favored, whatever $\hat{\mu}$.

In this context, the Bayesian perspective is correct and intuitive. 
This approach is useful in many contexts where we have a precisely defined alternative hypothesis.
However, this interpretation of the alternative hypothesis is \textit{not} what the authors intended. Although we wished  (i) \textit{to allow a large range of putative parameter values},  we also unintentionally specified a corollary: (ii) \textit{a vanishingly small prior density on the parameter manifold}.
In our conception of the statistical problem, 
we are \textit{not} interested in testing any \textit{precise model} for the distribution of $\mu$ (\textit{e.g.}~diffusion, stage drift, etc) as a  requisite for  determining whether the bead movement can be detected. 
If possible, we wish to achieve condition (i) without the corollary (ii). The predictive formulation of inference can achieve this goal. 
The vanishingly small prior density subtracts out of the predictive cross entropy as illustrated in Sec.~\ref{sec:dd}.




%

\subsection{Circumventing the Lindley paradox}
By partitioning the data into a training and generalization partition in the pseudo-Bayes factor, we are able to circumvent the most severe forms of the Lindley paradox by generating inference that is prior independent (for sufficiently large sample sizes). 
The postdictive limit depends sensitively on the data partition, but the predictive limit does not. Fig.~\ref{valvtrain} shows that for sufficiently large sample size $N$, the complexity rapidly converges to its limit as $\nu \rightarrow 0$ for $N_G<N_T$. Due to this convergence, two researchers  will report the same predictive pseudo-Bayes factor, even if they make different decisions about the prior and the data partition.

Our discussion of the Lindley paradox focusses mainly on critiques of a \textit{Bayesian} perspective and the defense of a \textit{frequentist} perspective on hypothesis testing or model selection.
In fact the \textit{frequentist} perspective we discuss includes methods from all three paradigms of inference.
Our criticism of the Bayesian paradigm is confined strictly to a criticisms of the use of Bayes factors and their undesirable consequences on model selection, as described above.
However, the Bayesian paradigm offers many strengths. 
The posterior is an elegant and intuitive framework for representing parameter uncertainty. Furthermore, hedging between parameter values (and models) typically leads to superior frequentist performance relative to point estimates.
Finally, the Bayesian paradigm offers a coherent framework for combining different types of data. Therefore we advocate retaining as many of these advantages as possible while eliminating paradoxical behavior in the context of model selection. The  pseudo-Bayes factor has  these desired properties.

Predictive methods (the information-based paradigm and predictive pseudo-Bayes factor) also circumvent many criticisms of the classical frequentist procedure: (i) Observed data that is unlikely in both the alternative and  null hypothesis results in the rejection of the null hypothesis.  (ii) An \textit{ad hoc} confidence level must be supplied. (iii) Only  pairwise comparisons between models can be made. (iv) A null hypothesis must be defined.  The predictive approach  circumvents each of these criticisms. Predictive methods also have provable asymptotic efficiency in terms of cross entropy loss in typical modeling situations\cite{Shao1997}, a feature which we  discuss in the supplement.


\idea{Loss of resolution}
To place the discussion of experimental resolution in context, it is useful to remember that biophysicists will routinely pay
thousands of dollars more for a 1.49 NA versus a 1.4 NA objective with nominally a 6\% increase in the signal-to-noise ratio. This obsession with  signal-to-noise ratio might suggest that a similar effort would be expended to optimize the resolution of the experimental analysis to exploit the data as efficiently as possible, especially in the context of single-molecule experiments where the sample size is often extremely limited. The Bayesian formulation of inference can imply a prohibitively stringent significance level for the discovery of new phenomena. The frequentist formulation of inference is tuned for discovery in the sense that it explicitly controls for the largest acceptable false positive probability. The Bayesian may require a much larger sample size to detect the same phenomena, as illustrated in Figs.~\ref{fig:Res}, \ref{Fig:Lindely} and \ref{fig:pval}.

\subsection{The multiple comparisons problem}
We have demonstrated that predictive inference has a lower threshold for discovery, but proponents have argued that the loss of resolution is in fact a feature rather than a flaw of the Bayesian paradigm. There is a perception that the canonical frequentist significance test is too weak and leads to spurious claims of discovery. An important and subtle problem with Frequentist significance testing is the multiple comparisons problem (multiplicity). For instance, if 20 independent false hypotheses for tumor genesis were independently tested at a 95\% confidence level, one would expect spurious support for one of these hypotheses. Multiplicity can arise in more subtle contexts: Hypotheses (or priors) are  modified after some results are known in the course of research, often unconsciously.  In singular statistical models, there is often implicit multiplicity in the maximum likelihood estimation procedure \cite{watanabe2009,CPlong,FICshort}. The peer-review process itself may favor the most extreme results among multiple competing articles. In exact analogy to the tumor genesis example, multiplicity can result in the spurious selection of the alternative hypothesis in each of these cases.

These false discoveries are a consequence of using an incorrect frequentist significance test  \cite{Dunnett1955}. For instance, we have described how the complexity in information-based inference must be modified in the context of a singular model \cite{CPlong,FICshort}. (\textit{e.g.}~\cite{Ioannidis:2005rr}). From a frequentist perspective, the significance test must reflect the presence of multiple alternative hypotheses which leads to  corrections (\textit{e.g.}~Bonferonni correction \cite{Dunnett1955}). These corrections increase the critical test statistic value to reflect the true confidence level of the test in the context of multiple alternative hypotheses.  In summary, the failure of frequentist methods due to un-corrected multiplicity is not a flaw in the frequentist paradigm but rather a flaw in its application. Bayesian inference can naturally circumvent some of these problems in a principled way,
but in many 
applications there are 
parameters for which one must supply an uninformative prior. As a result, the effective confidence level is \textit{ad hoc}. 
If multiplicity is the source of spurious false discoveries, a principled approach is to correct for this problem explicitly.

\subsection{Statistical significance does not imply scientific significance}

Simpler models  are often of greater scientific significance.  Therefore there is a perception that frequentism is flawed because it typically assigns higher statistical significance to larger models, relative to the Bayesian paradigm. This perception conflates statistical and scientific significance. Almost all natural systems appear to be described by models with a clear hierarchy of effect sizes \cite{Machta:2013hl}. Scientific progress is  achieved by studying  the largest effects first, irrespective of the statistical significance of smaller effects. 
The selection of effects to include in a model is a matter of judgment and scientific insight. There are important non-statistical systematic sources of error that must be considered. If sample size is large enough, these systematic effects will suggest the use of a larger model from a predictive standpoint, even if the larger model is not scientifically relevant \cite{marden2000hypothesis}.
Statistics  supplies only a lower bound on scientific significance by determining whether a hypothetical effect can be explained by chance.

\idea{Conclusion}
Bayesian inference can be powerful in many contexts, especially in singular models and in the small-sample-size limit where point estimates are unsuitable \cite{watanabe2009}. 
But Bayesian inference can result in strong conflict with frequentist inference when uninformative priors are used.
When a Bayesian analysis is desired, we advocate the pseudo-Bayes factor method \cite{gelfand1994} for inference on model identity  with a small ratio $\nu$ of generalization to training set sample size.
We demonstrate that only in this predictive limit can  inference be expected to be consistent between independent analyses. 
This approach is fully Bayesian for parameter inference, but free from the Lindley paradox. Therefore it preserves all of the advantages of Bayesian methods without the risk of paradoxical inference on model identity and optimizes experimental resolution.

\section*{Acknowledgments}
P.A.W. and C.H.L. acknowledge helpful discussions with S.~Presse, and M.~Drton, and constructive feedback from reviewers. C.H.L. would like to thank D. Mayo for describing the problem.
This work was support by NSF grants NSF-PHY-084845 and NSF-MCB-1151043-CAREER.

\bibliography{colinharrylamont}

\newpage

\appendix
\section{Calculation of Complexity}
In the large sample-size limit of a regular model, we can view a model as a flat prior on some subspace $\Theta$ of dimension $K$ embedded in a larger parameter space $J+K$. The marginal likelihood of $N$ measurements with dimension $J+K$ can then be written as;
\begin{align}
q(X^N) &= \integral{^J  \theta_\perp}{} \delta^J(\theta_\perp-\theta^0_\perp) \integral{^K \theta_\| }{}  \frac{\exp {\frac{-\sum_i(X^i-\theta)^2}{2 \sigma^2}}}{( 2 \pi \sigma^2)^{N(K+J)/2}}
\intertext{Which gives for the predictive distribution,}
h(X^{N_G}|X^T) &= h(X^N) - h(X^{N_T})\\
 &= \frac{\sum_i (X^i_\perp - \theta_\perp^0)^2}{2 \sigma^2 } + \frac{N_G(J+K)}{2} \log 2\pi \sigma^2 \\
 &\quad + \frac{K}{2} \log \frac{N}{N_T} + \frac{S_N - S_{N_T}}{2 \sigma^2} \nonumber
\end{align}
Where $S_{N_T}$, is the sum of (projected) squared deviations from the (projected) mean $\mu_{N_T} = N_T^{-1} \sum_{i\in N_T} X^i_\|$. A straightforward calculation shows that 
\begin{align}
S_N-S_{N_T} -S_{N_G} &= \frac{N_T N_G}{N} \paren{\hat{\mu}_{N_T}-{\hat{\mu}_{N_G}}}^2
\end{align}
Expanding around $\theta_0$, where $\theta_0$ is the parameter in the manifold $\Theta$ minimizing the KL divergence from the true distribution $p(\cdot)$,
\begin{align}
h(X^{N_G}|X^T) &= h(X^{N_G}|\theta_0) + \frac{K}{2}\log\frac{N}{N_T} - \frac{N_G^2}{N }\frac{\paren{\hat{\mu}_{N_G} -\theta_0}^2}{2 \sigma^2}\, +  \\
&\quad \frac{N_G N_T}{N} \frac{\paren{\hat{\mu}_{N_T} -\theta_0}^2 - 2\paren{\hat{\mu}_{N_G} -\theta_0}\paren{\hat{\mu}_{N_T} -\theta_0} }{2\sigma^2}.\nonumber
\end{align}
The deviance terms cancel under expectation. After rescaling, we can write,
\begin{align}
\frac{N}{N_G}\overline{h(X^{N_G}|X^T)} &= \overline{h(X^{N}|\theta_0)} + \frac{K}{2} { \frac{N}{N_G} \log\frac{N}{N_T}}.
\intertext{Defining $\nu = N_G/N_T$ emphasizes the limit behavior}
\frac{N}{N_G}\overline{h(X^{N_G}|X^T)} &= \overline{h(X^{N}|\theta_0)} + \frac{K}{2}(1+\nu^{-1}) \log \paren{1+\nu}.\label{eq:crossExpectation} 
\end{align}
The term $\overline{h(X^{N}|\theta_0)}$ is estimated by the observed information at MLE $h(x^{N}|\hat{\theta}_x)$. The error in this estimator (training error) is again $\frac{1}{2}\chi^2(K)$ distributed \cite{akaike1773}, making the following estimator unbiased
\begin{align}
\overline{h(X^{N}|\theta_0)} & \, \hat{ =} \, h(x^N|\hat{\theta}_x) + \frac{K}{2}. \label{eq:truthestimator}
\end{align}
When Eqn.~\ref{eq:truthestimator} is used with Eqn.~\ref{eq:crossExpectation}, it gives us an information criterion corresponding to the pseudo-Bayes Factor for each partition choice $\nu$.

\section{Definitions and Calculations}\label{app:volume}
\subsection{Volume of a distribution}
Our intuitions about the volume of a distribution can be made mathematically precise using the self-entropy $S$. The self entropy functional is defined as
\begin{equation}
S\left[q(\cdot)\right] \equiv - \integral{{\bf \theta}}{} q({\bf \theta}) \log q({\bf \theta}),
\end{equation}
and the volume is defined in turn as
\begin{equation}
V_q \equiv e^{ - S\left[ q(\cdot)\right] }.
\end{equation}
For uniform distributions, this entropic definition reduces to the volume of the support. 
For a normal distribution of dimension K the volume is
\begin{equation}
V_{\Sigma} = \paren{2 \pi e}^{\frac{K}{2}} |\Sigma|^{\frac{1}{2}} \approx \paren{2 \pi e \sigma^{2}}^{\frac{K}{2}}.
\end{equation}
where the second equality only holds if $\Sigma$ is proportional to the identity.

\subsection{Showing $I = I'$ to zeroth order}

The first term in the information difference
\begin{align}
I'(x^N)&= \hat{H}^{1|N-1}-\hat{H}^{N|0}\\
 &= \paren{ \expectation{\theta}{\varpi(\cdot|x^N)}  h(x^N|\theta) } - h(X^N) + {\cal O}(N^0)\\
&= \expectation{\theta}{\varpi(\cdot|x^N)}  \log \frac{q(x^N|\theta) } {q(x^N)} + {\cal O}(N^0) .
\intertext{By multiplying the numerator and denominator by $\varpi({\bm \theta})$, we can identify this first term as the KL divergence that we used to define $I(x^N)$}
&= I(x^N) + {\cal O}(N^0) 
\end{align}

\subsection{Significance level implied by a data partition}

Under the assumption of the smaller model, the information difference is expected to be distributed like a $\frac{1}{2}\chi^2$ with $\Delta K$ degrees of freedom, where $\Delta K = K_2-K_1$. The effective significance level $\alpha_\nu$ is therefore
\begin{equation}
\alpha_\nu = 1 - \mathrm{CDF}\left[\chi^2_{\Delta K }\right]\left(\Delta K \left[1+(1+\nu^{-1}) \log (1 + \nu)\right]\right).
\end{equation}
This function is plotted in Fig.~\ref{fig:pval} for two choices of the dimensional difference. An interesting corollary is that for large $\Delta K$, typical confidence levels may actually be less than equivalent predictive methods such as AIC.  In other words, we can reject the null hypothesis before it become predictively optimal to use the larger model.

\section{Other Methods}\label{app:othermethods}
There are several methods that deviate more drastically from the standard Bayesian probability calculus. We mention here just a few of the interesting ideas which have been proposed.

\subsection{Other Predictive Estimators}
Once a data division strategy is chosen and we can agree on what we are trying to estimate, there are many information criteria which can be used. 
For instance, the predictive limit can be estimated using AIC, DIC \cite{Spiegelhalter2002} and WAIC \cite{watanabe2009}. When the posterior is known to be approximately normal, AIC can perform with minimal variance \cite{Shibata1997}. 
Far from normality and the large sample-size limit, WAIC has a uniquely well developed justification in terms of algebraic geometry, but the standard LOOCV seems to have better properties in numerical experiments \cite{Vehtari2015pareto}.
Similar alternatives to BIC exist for postdictive performance estimation \cite{Leung2014,watanabe2009}.

\subsection{Data-Validated Posterior and Double use of Data}
Aitkin \cite{Aitkin1991} attempted to address the Lindley paradox by proposing training and validating the data using the entire dataset $X^N$. The resulting posterior Bayes factor comes from the observed posterior information:
\begin{equation}
H_\mathrm{POBF}(X^N) = h(X^N|X^N)
\end{equation}
This has a complexity ${\cal K}_\mathrm{Aitkin} = \frac{K}{2} \log 2 \approx 0.35 K$. This is far too weak to realize Occam's razor. This weakness results from two effects: i.) We use here the a generalization sample size of $N_G = N$ instead of the predictive limit where the generalization sample size is zero. ii.) The double use of data means that the posterior is over-fit to the particular dataset. This posterior appears to performs better than even knowledge of the true parameter ${\cal K}_0 = \frac{K}{2}$. 
Overfitting can also occur when data are double used implicitly through posterior training, as in empirical Bayes methods where prior hyperparameters are optimized with respect to the model information.

We do not believe that the double use of data is completely ruled out of a principled statistical analysis \cite{akaike1978}.
But because double use of data is antithetical to the interpretation of conditional probability, and because it very often leads to overfitting, double use of data requires careful justification.

\subsection{Fractional Bayes Factor}
O'Hagan \cite{ohagan1995} sought to define a minimal training set mathematically by taking some small power of the likelihood. The fractional model information is then
\begin{equation}
H_{\mathrm{FBF}}(b) =\log \expectation{\param}{\prior} q^{b}(X^{N}|\param) -\log \expectation{\param}{\prior} q(X^{N}|\param)
\end{equation}
where $b$ is chosen to be very small. If epsilon goes to zero, this expression is obviously identical to the original model information. 
As O'Hagan notes ``The key question remaining in the use of FBFs is the choice of $b$. It may seem that the only achievement of this paper is to replace an arbitrary ratio [{\em i.e.} $N_G / N_T$] with an arbitrary choice of $b$." The same issues with defining a minimal experiment for minimal training also arise for this approach.

\section{Efficiency and correct models}
The landmark treatment by J. Shao\cite{Shao1997} and its discussion by Yang \cite{Yang2005} are sometimes viewed as supporting BIC and Bayes factors in certain situations. We therefore wish to discuss this important work in more detail. We suppress many of the technical details for the purposes of our discussion, and refer to \cite{Shao1997,Yang2005} for more precision.

Let $\alpha_0^N$ be the identifier for the most predictive model at sample size $N$ (which may not be the true model!), and let $\hat{\alpha}_\nu(x^N)$ identify the model chosen by selecting the largest pseudo-Bayes factor parameterized by $\nu$. We can define the loss ratio in terms of the predictive cross-entropy of the trained model,
\begin{align}
\epsilon_\nu(x^N) \equiv \frac
{\expectation{Y}{}
	h(Y |x^N, \hat{\alpha}_\nu(x^N) )}
{\expectation{Y}{} h(Y |x^N, \alpha^N_0) },
\end{align}
where the expectations are taken with respect to the true distribution. Shao identifies a reasonable criteria for the performance of a model estimator $\hat{\alpha}$: {\em asymptotic-loss-efficiency} which is equivalent to the condition that
\begin{align} \epsilon_\nu(X^N) \rightarrow_p 1,
\end{align}
That is, the loss ratio converges in probability to unity as the sample size goes to infinity. 

Shao found that the {\em context} in which model selection is performed is incredibly important to whether or not asymptotic efficiency is achieved. Specifically, there are two very different situations:
\begin{enumerate}
\item There is no correct model, or there is exactly one correct model which is not nested inside a more complicated model.
\item There is more than one correct model. The smallest correct model is nested inside potentially an infinite set of increasingly complicated models, which are all capable of realizing the smaller model.
\end{enumerate}
If condition (1) holds predictive methods ($\nu\rightarrow0$) are guaranteed to be asymptotically efficient, and (pseudo-)Bayes factors for which $\nu>0$ are {\em not} guaranteed to be asymptotically efficient. But if condition (2) holds, then statistical fluctuations will cause AIC and pseudo-Bayesian methods to choose larger models than $\alpha_0$ with a probability that never goes to zero. It is necessary for the penalty that is, $\nu$, to diverge to ensure that the probability of choosing a larger correct model will converge to zero, and that asymptotic efficiency can be achieved.

If the possibility of condition (2), and the true model is realizable at finite dimension, many would suggest that we are justified in using Bayesian methods which have a divergent $\log N$ penalty and thus hope for asymptotic efficiency. We criticize this position on several points.

First, condition (2) is unlikely to ever hold. The Boxian proverb, ``All models are wrong,'' expresses the general truth that nature is too complicated to ever yield the exact truth to a finite dimensional model. Condition (1) is far more likely in any typical modeling context.

Second, whereas predictive methods occupy a unique place in relation to condition (1), the rate at which penalties must go to infinity to satisfy efficiency under condition (2) is not uniquely determined. All methods whose complexities go to infinity slower than $N$, will (with some technical caveats) satisfy asymptotic efficiency. A complexity of $\log \log N$ would be no less favored under this argument than the $\log N$ complexity of BIC.

Finally, the asymptotic argument which prefers BIC under condition (2) seems to have little bearing on the conditions we would observe at finite sample size. At finite sample size, we do not know if we are in the regime where we are selecting the true model, or if the true model cannot yet be resolved with the available data. If the true model cannot be resolved, we'd still expect AIC to typically outperform BIC for the same reasons that hold in condition (1). BIC is unjustified unless we know a priori the scale at which a true effect will be observed. This is exactly the situation which holds when we have a precise distribution for the parameter of interest, and the Bayesian approach is indistinguishable from the way the frequentist would use {\it a priori} information in accordance with the Bayes law.


\end{document}